\documentclass[a4paper,12pt,hyphens]{article}

\usepackage[backend=biber,style=authoryear,sorting=none]{biblatex}
\usepackage{csquotes} 
\usepackage{hyperref}

\usepackage[hyphenbreaks]{breakurl}

\usepackage{amsmath,amsthm,amssymb}
\usepackage[all]{nowidow}
\nocite{*}

\newtheorem{thm}{Theorem}

\newtheorem{prop}{Proposition}
\newtheorem{cor}{Corollary}

\newcommand\nl{\hfill\break}

\newcommand{\eop}{\hfill $\square$ \vskip10pt}
\newcommand\ora[1]{\overrightarrow{\text #1}}
\newcommand{\an}{$\angle$}
\newcommand{\cg}{$\cong\ $}

\newcommand{\pr}{${}^\prime$}
\newcommand{\st}{${}^*$}
\newcommand{\dg}{${}^\circ$}

\begin{document}

\title{An analysis of Euclid's geometrical foundations} 
\author{Peter Malcolm Johnson\thanks{Universidade Federal de Pernambuco, Recife, Brazil. Email: \href{mailto:peterj@dmat.ufpe.br}{peterj@dmat.ufpe.br}}}

\maketitle

\begin{abstract}
The initial techniques developed in Euclid's {\em Elements}, well  before the use of the parallel postulate,
are reexamined in order to clarify even the most obscure details, particularly those related to equality,
superposition and angle comparison.  Some commentary on modern developments is included. The known but 
often misunderstood implicit handling of betweenness and points of intersection is briefly treated.
We also sketch a rigorous treatment of absolute geometry in a spirit similar to Euclid's, one that 
allows properties of angles and triangles to be derived from two simple axioms on right angles, 
which then leads to rigid motions of certain planar geometries.
\end{abstract}


\


\parskip5pt

\section{Some remarks on Euclidean geometry} 

Our ideas that arose from reading parts of the {\em Elements} and related works
are those of a mathematician with some knowledge of the historical background,
not those of a historian of mathematics familiar with the subtleties of Greek language and thought.
In what follows, there will be a fruitful interplay between different viewpoints, ancient and modern, textual and mathematical.
The principal motivation is to understand Euclid's approach at the beginning of the {\em Elements}, where the first foundations
of plane geometry are developed.   We provide and defend a critical  interpretation in which defects are examined but no obscurities remain, 
thus making this part of Euclid's text fully comprehensible by a large class of readers.  
Some suggestions for alternative approaches are made, especially towards the end.
Other commentary, often terse, is interspersed in an attempt to situate the present work 
within the large body of related research, on both ancient and modern themes. 

The {\em Elements} of Euclid, despite its limited scope, has had a deep and long-lasting influence on the development of mathematics 
and also on the history of mathematics.
The meanings and purposes of its content, and of other surviving versions and fragments of Greek and Hellenistic mathematical writings,
have long been the object of incisive analyses and debates, one important theme being the role of construction.
Although modern scholarship has advanced, Heath's translation of Euclid [Eu],
based on Heiberg's Greek edition,
is accompanied by a wealth of information that synthesized all relevant knowledge at the time.
The commentary on Book I, written by Proclus [Pr] about seven centuries after Euclid, contains much of interest.
It is a uniquely detailed study of Euclidean geometry from a philosophical perspective.  
We use Morrow's edited translation of that work but, as is customary, page references are to Friedlein's Greek edition.
 
In Euclid's presentation of plane geometry, 
conclusions are drawn from the study of idealized versions of geometric diagrams formed from finitely many points (roughly, 
markers of position, often where objects cut or touch each other), straight lines (segments), and parts of circles.  
Angles, another important ingredient, will be treated later. 
There is an implicit notion of incidence, expressed in various ways, but not of course the conception of objects as sets of points, 
used below only for ease of expression.
To prove results, suitable new points can be selected (they come into existence on the diagram), and straightedge-and-compass
constructions made, adding to the diagram but erasing nothing.
Such arguments are completely constructive, except those that make use of the Fifth Postulate on parallels.
Here we diverge from Bl\aa sj\"o's constructive thesis [Bl], which makes no such exception.
 
Euclid's approach is not static, as geometric procedures can be carried out.
To maintain the philosophical view of an ideal, unchanging geometric world, Proclus [Pr, 78] 
explains that all such procedures are acts taking place in the imagination.
Applications are another matter:  ``... when it touches on the material world ..." [Pr, 63].
The prevailing view is that a few proofs in the {\em Elements} rely on moving parts of diagrams.
That interpretation is questionable, as it is sufficient to construct copies of objects in specified positions.
More on this much-debated topic will appear below.

Relatively recent approaches to axiomatic Euclidean geometry usually work with some form of Tarski's system.
Among the highlights are Gupta's extremely ingenious constructions of various geometric objects, including 
midpoints of segments. This material forms part of the influential book [SST].
Approaches to constructive geometry in similar systems remain an active area of research.
Some of these have been formalized in appropriate languages that permit computer-verified proofs.
The article [BNW] gives an excellent exposition of how a meticulous analysis of a system like Euclid's is needed 
in order to be able to carry out the complex project of complete formalization.  Our much less formal approach 
to the most basic results of absolute geometry has key differences, notably concerning axioms for right angles, 
which may be of interest for future explorations in the area.  The main purpose, however, is to propose a way
to fully clarify the initial content of the {\em Elements}, suggesting modifications where necessary to ensure rigor.  

We now focus on how the most basic machinery is developed early in Book~1, how it can be further exploited in new ways,
and what assumptions, especially implicit ones, are involved.  
These analyses will eventually lead to a way to obtain a distinctly different alternative approach to axiomatization and exposition,
rigorous but still faithful to the spirit of Euclid's geometry, one that also has notable consequences of a modern form.

\section{Lines and their points}

To avoid confusion, what Euclid calls a straight line will often be called a {\em segment},
as segments of circles will not be discussed.
A segment AB is defined uniquely by specifying two distinct points, its {\em ends} A and B.  
As a line, it is extendible at either end, and  retains its identity, in some sense, if the ends change.  
Lines are never diminished, presumably because parts of diagrams cannot be erased.
By explicit assumption, non-extreme points can be chosen on any segment AB.
Such points are precisely those said to lie between A and B, in the strict sense of betweenness.  
It will soon be seen how this relation can be clearly understood in what we believe is its proper context.
  
Contrary to what one might expect, Proclus [Pr, 285] seems reluctant to admit that there could be
an idealized notion of an infinite line, extending beyond any given bound.
He argues that this cannot be the object of knowing imagination. 
In Euclid, an infinite straight line is treated as an arbitrarily extendible segment with unspecified ends.  
Although not clearly stated, any two such lines are supposed to have at most one point in common 
if they differ, a view repeatedly defended by Proclus against objections that had been raised by Zeno and others.
A convenient compromise, adopted below, is to say that two segments {\em overlap} if they share more than one point,
while segments are {\em collinear} if some other segment overlaps each of them. 

The Fifth Postulate avoids mentioning infinite lines, which rarely appear in the {\em Elements}. 
It often suffices to extend a segment using a bound that can be determined from the diagram.
Euclid's Prop.~I.27 gives an angular criterion showing that any point P not on a line $l$ lies on a line parallel to $l$.
The so-called Playfair form of the parallel postulate asserts the unicity of this (infinite) line, whereas
Euclid's form, the Fifth Postulate, asserts that any other line through P intersects $l$, in the direction expected from the sum of two angles, 
provided the two lines (as segments) have been sufficiently extended.
Our opinion, one among many, is that this postulate caused much discomfort because of the lack of 
a constructive procedure, not even one with an indefinite number of steps,  for obtaining the point of intersection.
We remark that in a non-archimedean geometry one could repeatedly double the two line segments at the correct ends 
(replacing any AB with the segment AC having midpoint B) without necessarily obtaining a point of intersection.
The constructibility of intersection points is a subtle matter that depends on the exact choice of axioms adopted and even on the underlying logic.
Boutry [Bo], in Chapter I.3, gives a detailed analysis, and among much else uses metamathematical methods to give a new proof of 
the independence of the parallel postulate in Tarski's system.  That postulate can in fact  be formulated in many different ways  
([Bo, App. D] lists 34 of them), and their logical dependencies depend on which axiomatic framework is used.

It is now well known, although perhaps not among scholars of ancient mathematics, that all expected results about betweenness, or of 
the ordering of points on a given line, follow easily from a ``plane-separating" assumption often made in Euclid's proofs about plane geometry.
The method was discovered a few years after Hilbert's book [Hi] appeared in 1899.
The history of related developments at that time is ably recounted in [BB].     
Some form of Euclid's assumption (often that of Pasch) is adopted as an axiom in contemporary geometry textbooks that treat foundations. 
The original assumption, made explicit and recast in modern language, is that
any infinite line $l$ determines a partition of the points not on $l$ into two classes, called half-planes or the sides of $l$, 
and each side is convex: it contains all points of a segment AB whenever it contains A and B.

The key observation is that, given a point O on a line $l$, the sides of any other line $m$ that passes through O will separate $l$
into two rays (half-lines) based at O and otherwise disjoint.  
Moreover, this separation does not depend on the choice of $m$, as any two points A, B on $l$ other than O 
lie on different sides of $m$ precisely when O lies on the segment AB.  Proclus [Pr, 198], referring to a work of Pappus, 
comments that ``... a point divides a line, a line divides a plane, and a plane divides a solid ...",  just after an assertion 
about the homogeneity of lines and of the plane, all of which supposedly follow from axioms and definitions.

A closely related assumption, used implicitly by Euclid, is that a segment with ends on different sides of a line cannot fail to intersect that line.
As indicated in some textbooks such as [Gr], [Ha] or [Mo], the simplest form of Pasch's axiom, about how lines and triangles can intersect,
is equivalent to Euclid's assumptions on plane separation.
One can continue the study of O on $l$ by letting O vary over the points of $l$.
This leads easily to a definition of two mutually opposite total orders on the points of $l$, and each of these orders determines the segments of $l$. 

Modern geometry texts that contain material related to order often follow a sub-optimal didactic sequence, 
largely due to the influence of Hilbert's book [Hi].
An accessible non-metric approach is to impose axioms about binary order relations (paired with their opposites)  for the points on each line, 
then use the axioms to define segments and to derive properties of betweenness.

Although the so-called continuity assumptions about intersections have been criticized for centuries by readers of Euclid, De Risi [DR1]  
presents a spirited and well-founded defense of such Euclidean practices. To give a crude paraphrase of a sophisticated argument, 
points can be created on diagrams to mark positions where parts of figures touch or cross.  
The view that they must be selected from a pre-existing background of points is a modern one.
We remark that the interior of any triangle is convex, as it is the intersection of three half-planes. 
The intersection of this interior with a line, if non-trivial, must then form a bounded interval. 
It is a natural step to use, or create where necessary, points to mark the ends of that interval where the line crosses the triangle that 
bounds the interior.  Similar remarks apply to interiors of circles, once their convexity has been proven.

\section{Magnitude, congruence and superposition}

The assumptions that geometric objects have (or even are) magnitudes, and that those of the same kind can be compared, 
go back to the very roots of geometric thought. 
Results about associated concepts such as ratio and proportion are developed in the {\em Elements}, but in that work
the fundamental concept of magnitude is taken for granted, notably in the common notions, rather than made explicit.
We do not know what attempts can be found, in the vast literature on the subject, to clarify this interesting historical detail.  
Either explanatory comments were never included or they were omitted at an early stage,
presumably because such ways of thinking were so commonplace at that time.
Writings of Aristotle, and the later commentaries of Proclus, cast light on these and other facets of educated Hellenistic thought.
 The influence on ancient geometry of the acceptance or not of objects with infinitesimal magnitude (not just horn angles)  
is a topic we do not enter into.   

In modern language, congruence, for geometric objects such as segments or angles, is required to be an equivalence relation $\cong$
whose classes support a total order relation, subject to some natural requirements. 
The order pulls back to a pre-order $\leq$ on the objects themselves, with congruence in place of equality.  
The associated notation $>$ is also used.
 The relation called equality in Euclid coincides with this notion of congruence, at least for segments, angles and, in certain contexts, triangles.
Equality for regions, usually bounded by polygons, is a weaker relation than congruence and will not be discussed here.
In all cases, ``equality"  seems to mean equality of an appropriate ``substance" or geometric magnitude, 
something whose definition lay beyond the expressive capacity of Greek numbers.
It remains to be seen how usages in the {\em Elements} point to implicit mathematical definitions for relations of equality and order,
where those definitions will depend on the type of geometric objects being compared.

 The common notions (we follow Heath's use of Heiberg's reduced list of five) encapsulate basic properties that magnitudes are supposed to satisfy.
They, along with definitions of  geometric  objects such as points, lines and angles, form part of the preliminary material of the 
{\em Elements} whose original form, if indeed there was a single original, is unknown.
De Risi [DR2] provides an extensive survey of the literature on the common notions, showing what is known in a wider historical context.
He presents abundant textual evidence to support the long-held view that some of these notions
are later additions, as they are not always applied, or cited, in expected ways.  
The assumptions separated out as common notions do, however, play an essential role in the {\em Elements}.
While propositions are cited indirectly in proofs by inferring their conclusions immediately after verifying their hypotheses, no uniform form 
of citation for common notions can be expected, as the precise details will depend on the relevant parts of a diagram under consideration.
We will closely examine how these notions are, or could be, used in proofs of the initial and most fundamental results in the {\em Elements}.

We first focus on the Fourth Common Notion, whose English translation as ``Things which coincide with another are equal" needs clarification.
Since definitions of equality cannot be found elsewhere, we 
believe that this is not merely a criterion for equality,
but the very definition itself, as will be justified below.  The Greek verb used here is the same that appears, inflected differently, 
in three places in the {\em Elements} (Propositions I.4, I.8 and III.24), where it is usually read as an action of superposition.
Heath [Eu, p.~225] states that Euclid's phraseology leaves no room for doubt that a figure is moved and placed on another.
Sidoli [Si, Sec.~5.3], argues that the verb is used in the genitive absolute and only implies a hypothetical situation.  To paraphrase,
it is not that a part of a figure is moved to coincide with another, but that it could be.  
The article of Axworthy [Ax] shows the intensity,  several centuries ago, of conflicting arguments about this topic. 
Deferring until later a critical analysis of the cumbersome process of transport used in the {\em Elements}
for achieving superposition, we now give a brief overview.

Given a segment BC and a point A, the aim of Euclid's Proposition I.2  is to construct a segment equal (congruent) to BC with an end at A,
in effect a copy of BC.  The copying procedure is often called transportation in what follows, although the original segment is not erased
and it seems clear that the purpose is to avoid any need for actual motion.  
As its first application, in I.3, different segments are compared and the smaller is cut off (subtracted) from the larger.
For practical purposes, the construction legitimizes motion by simulating it, and also allows segments to be compared.
In particular, it provides the definition for equality (congruence)  of segments.  
Its fundamental importance can be inferred from its very early placement in the {\em Elements} and the 
considerable effort needed to effect it, one that combines the use of a triangle (obtained in I.1 from two circles), 
extensions of segments, two more circles, and some common notions about operations on segments.
In constrast, contemporary textbooks treat congruence and the associated common notions, for segments and angles, 
as parts of the axiom schemes adopted, so that transportability becomes irrelevant.

A trivial but important consequence of I.2 is that a triangle can be copied, or superposed, onto another (or the same) triangle
when the corresponding sides are congruent.  Euclid provided a form of the expected SSS construction in I.22, citing I.20 
(the triangle inequality) to give necessary conditions on side magnitudes, assuming implicitly that certain circles intersect. 
This is followed immediately by I.23, which provides the SAS criterion for copying a given angle (really as part of a triangle) to a 
specified location. It appears from this and the very early I.4, often used later, that Euclid placed special importance on angle transport, 
a method that will  be explained later.

For our foundational purposes, it suffices to impose a weak SSS axiom, just below, to ensure the transportability of an existing triangle to any 
location, as specified precisely in the axiom.  In the {\em Elements}, a construction of this provisional form would fit well immediately after I.2. 
Although that was not done, we argue later that the idea, used implicitly, is what justifies arguments involving superposition.   
(In III.24, use a triangle formed from a chord and a circle center.)
The usual SSS criterion for the congruence of triangles cannot be adopted until I.7, the crucial uniqueness or``rigidity" result, has been established.

\medskip
{\bf Axiom T:  If A, B, C are non-collinear points and A\pr B\pr\ \cg AB, where  A\pr\ and  B\pr\ lie on a line $l$, 
      then on either side of $l$ there is a point C\pr\ with  A\pr C\pr\ \cg AC and  B\pr C\pr\ \cg BC.}

\section{Equality and comparison of segments}

We return to our critical study of Euclid's foundations by considering how the relation of congruence (\cg) for segments could be formalized axiomatically.   
If it is desired to do this without (or before) imposing an axiom about how lines separate the plane, it suffices to require that distinct segments 
of the form OP and OQ, each overlapping (as defined above) with a fixed ``reference" segment OA, cannot be congruent.  
Comparisons of their magnitudes can be defined as follows.  
If Q lies on OP then OP$\,>\,$OQ (``the whole is greater than the part"). To ensure the converse, an axiom must be imposed.  
In modern language, where segments are regarded as certain sets of points, this amounts to saying that the set of segments OX 
that overlap with OA is totally ordered and that this relation coincides with that of set inclusion. The definition is clearly independent
of the initial choice of one OA among the segments OX.  The ends of these segments OX form the ray or half-line $\ora{OA}$ based at O. 
In addition, varying O, each point on a line divides that line into two opposite rays with a single point in common.  

Next, segments OP and OQ with a common end O are defined to be congruent precisely when some circle with center O contains P and Q.
It is assumed, implicitly, that every ray based at O intersects the circle at a unique point.  One thereby obtains a congruence-preserving 
order-isomorphism between the sets of those segments of $\ora{OP}$ and of $\ora{OQ}$ of the form OX.

The proof of Euclid's I.2 is one of several places where an equilateral triangle (ABD in the translated text) is used, where D is found by 
intersecting two circles.  Given that such a point D exists, it could be expected that there is another such point D\pr\ on the other side of AB.  
If  one is working with axioms that are sufficiently strong to allow proofs, at this very early stage, of the usual results about reflections in a line, 
the choice of D or D\pr\ in the construction will produce an equal (congruent) copy of the segment. 
A more serious flaw, despite Heath's attempt to minimize it [Eu, pp.~243 and 261], is that it has not yet been proven that 
different circles cannot have even more common points. 
The uniqueness of the copying procedure, an essential ingredient for developing further machinery needed in later proofs,
can be forced by adjoining an appeal to the common notion that the whole is greater than the part.
It is unsatisfying to have to make such a strong assumption to make this fundamental construction work in an unproblematic way.

We suggest a slightly different way to proceed.  To avoid constructing a triangle ABD, D can be taken instead to be the midpoint of AB.
One could posit the existence of midpoints of segments, axiomatically, just as lines and circles are supposed to exist.  
Issues about constructibility, and ways this could be proven later, are not of immediate concern. 

Euclid's method for copying a segment BC to A, expressed informally in language more dynamic than the text in I.2 supports, is to
rotate (if necessary) BC around B to the segment BE for which $\ora{BE}$ and $\ora{BD}$  are opposite rays.  
In other words, it is assumed that DB can be prolonged sufficiently at B to a segment DE, where E lies on the circle with center B and radius BC.
We do not see any ambiguity or need for case division, read into the text here by Proclus and some later commentators. 
Now DE can be rotated about D to a segment DF, taking A to B. Finally, AF is congruent to BC by the common notion on subtracting equals from equals.  
If instead the objective was to transport BC to a fixed ray based at A, another rotation, of AF about A, may be needed. 

This procedure, in essence, transports circles centered at one point B to circles of the same radius centered at another point  A.
It could best be interpreted as giving a definition of what it means for circles with different centers to have the same radius, 
or for segments to be congruent, and can then be used to compare magnitudes.  As the steps of the construction are reversible, 
the relation between segments is symmetric. The construction can also be applied to AB to interchange A and B.
Whether or not one wishes to invoke the Fourth Common Notion here, to something that coincides with itself, 
there is no doubt that magnitude for undirected segments should not depend on the order in which the two ends are specified.

Axioms are needed to ensure that segment congruence is transitive.  There may be several good choices, and we leave this open.
Euclid's solution was to simply rely on common notions about equality.  Proof of transitivity requires showing that a circle at some point 
A (its center) remains unchanged if it is transported to another point B, then to C, then directly back to A. 
   
With little effort one obtains the following strengthening of I.2, which we prefer to express 
in modern language and with a different naming convention.
What we find surprising is that Euclid did not see fit to include here some explicit result about  
the transport of segments that have been subdivided or extended, one that would justify segment arithmetic. 
   
\begin{thm}
The transport of a segment AB on a line $l$ to some A\pr B\pr\ on $l'$, taking A to A\pr\ and B to B\pr, induces a unique
congruence-preserving bijection $X \to X'$ between the points of $l$ and those of $l'$.
\end{thm}

{\em Proof:}  
 The construction in I.2, applied to both AB and AC, transports any segment AC that overlaps AB 
to the segment A\pr C\pr\ that overlaps A\pr B\pr\ and is congruent to AC.   
When C lies between A and B, the magnitudes satisfy  C\pr B\pr\ = A\pr B\pr\ -- A\pr C\pr\ = AB -- AC = CB.
A point D on the ray opposite to $\ora{AB}$ should be transported to the point D\pr\ on $l'$ with A\pr D\pr\ \cg AD
such that D\pr B\pr\ = D\pr A\pr + A\pr B\pr, giving  D\pr B\pr\ \cg DB.  Other cases have similar proofs. \eop
\smallskip

The common notions can now be exploited to define algebraic structures on lines.  
We make only some terse remarks to supplement mathematical developments like those in [Ha, Ch.~4].
Ignoring the long history of resistance to related ideas, it is advantageous to make the leap to ordered lines (not just rays) whose points 
form an abelian group of signed magnitudes.  Dependence on the choice of a zero in the group of points can be avoided by passing,
via a canonical isomorphism, to the group of congruence-preserving translations, which acts regularly on the line.
Reflections, which fix one point and reverse the order, are not used here.  The transition from one ordered line to another induces 
a unique isomorphism between the groups of translations.  If a fixed ordered line, with its associated ordered abelian group $T$, is chosen, 
each unordered segment can be assigned a certain positive element of $T$ as its magnitude, or as the distance between the points at its ends.    
Then segments will have equal magnitudes precisely when they are congruent.

When continuing with the analysis of $T$ to obtain an ordered field, a segment deemed to have unit magnitude is often chosen.
One can instead follow Nagumo's approach to positive quantities, translated in [Na].  
Here the operations of positive multiplication are by definition the order-preserving automorphisms of $T$.
This can of course be related to classical ideas on ratio and proportion.

\section{Angles in Euclidean geometry}
 
In a modern conception, an angle consists of two non-opposite rays, its sides, that are based at the same point, its vertex.  
At intersections of pairs of lines, each line $l$  cuts the other into two rays on opposite sides of $l$, by the convexity of plane-separation.  
The rays form four angles at their common vertex and each angle has a well-defined interior, determined by which side of each of the 
two lines contains it.  Inclusion of interiors defines a natural total order on the angles \an AOX that share a common ray $\ora{OA}$ and lie on 
the same side of the line determined by that ray.  There is also a partially defined operation of sum (and difference) for angles with 
a common vertex, where  \an AOB +  \an BOC  $=$ \an AOC$\,$  whenever B lies in the interior of \an AOC. 

To address the important topic of angles in the {\em Elements} and to pinpoint possible flaws in the treatment,
we make our own analysis. Critical analyses of this particular topic seem to be rare.
Wagner [Wa] is close to our point of view, Shabel [Sh, $\mathsection$1.3] partly so,
whereas Alvarez [Al] (the most detailed), Mueller [Mu, pp.~19--25] and Panza [Pa, pp.~91--99] are not.
That of [BNW] is very different, mainly because it avoids using the plane-separation assumption.
Sidoli [Si, Sec.~5.3] is useful for clarifying some subtleties of Euclid's uses of technical language.
 Only later will we develop an alternate and provably reliable method for treating angles.

Definition 8 in the {\em Elements}, of an angle as the mutual inclination of two lines, is unsatisfactorily vague, and much has been written about it.
While several kinds of angles were accepted and studied by geometers, only rectilinear ones (Definition 9) will be considered below.
They will be treated as mathematical objects, not for example as qualities or relations between lines, so important to Proclus but not to us.
A careful study of the ways angles are used in proofs is sufficient for allowing one to infer what they are, how they can be compared, 
and what properties they are assumed to have.

In the {\em Elements}, a named angle, say ABC, can be defined to be an object 
that consists of two non-collinear segments AB and BC, its sides, and has B as its vertex.  
This angle \an ABC coincides with \an CBA and forms part of a triple of angles in a triangle ABC.  
An angle presumably remains ``equal"  if its sides are extended, but it is undesirable to assume this uncritically,
as Euclid did in I.23 and in earlier places, mentioned below. 
In Tarski's system, the Five-Segment axiom handles this problem.  We will find another way.

To compare any two angles, the only plausible method is to copy angles as parts of triangles.  
Not until the proof of I.8 does Euclid reveal, by making explicit use of it, that this is his otherwise elusive definition for the 
equality of angles. It is  merely an application of the common notion about things that coincide (by superposition) with each other.
The uniqueness result I.7, on the rigidity of triangles, is needed to make this definition useful.
Any angle, as part of a triangle,  can then be copied uniquely to a given location. This, without uniqueness, is the content of I.23,
whose statement fails to mention that the construction relies on the prior choice of a suitable triangle.

The central problem is how to establish I.7 without assuming anything not yet proven.  A counterexample would consist of two distinct 
triangles,  copies of each other, lying above a common base.  For reasons that are not clear, it is not observed that the corresponding 
angles at the base would then be equal, forcing one in each pair of equal angles to be a proper part of the other, an absurdity.   
The actual proof of I.7 involves two cases, one handled in the {\em Elements} and the other by Proclus.
To expose and clarify the underlying assumptions, the proofs  of  I.4 and I.5 must also be examined.

Euclid's Proposition I.4 is the first theorem of the {\em Elements} after three initial constructions.   Given triangles ABC and DEF with 
AB \cg DE, AC \cg DF and equal angles at A and D, the objective is to prove that corresponding sides and angles are equal.  
The proof has been heavily criticized for its supposed illegitimate use of superposition and an unconvincing conclusion, 
as well as its failure to explain equality of angles.  We furnish evidence just below 
that the proof of I.4, as usually read even today, has been seriously misunderstood and is a caricature of the intended one.
It is noteworthy that Proclus, who devoted so much of his commentary to raising and dealing with all sorts of possible objections,
described and even praised the proof procedure.  Our reading of it is compatible with his description.

Under the above hypotheses of I.4, one can certainly copy AB to DE and, to ensure equal angles, copy AC to DF, 
but then nothing supports the conclusion that BC \cg EF.
A legitimate procedure is to apply (copy or transport) the triangle ABC,  not necessarily uniquely, to some DEG 
with F and G on the same side of DE.   To be able to exploit the hypothesis about equal angles, we believe it must be read 
as the assumption that some copy DEG of ABC has the same angle at D as DEF. 
Since DG \cg AC \cg DF, G is forced to coincide with F, and ABC (or a copy) has indeed been superposed on DEF.   
Each remaining pair of corresponding angles of ABC and DEF is then equal by definition, from angle transport.  

In order for I.4 to be usefully applied before uniqueness results have been established, it should be ensured that an angle can never 
be equal to some other in its interior, since Euclid's sequence of basic initial results depends on repeated use of this inference.
We believe it is best justified as an implicit invocation of the Fifth Common Notion, that the whole is greater than the part.
Mueller [Mu, p.~35] marshals evidence suggesting that the inference, as originally made, was instead a direct and intuitive one, 
made obvious by the diagram. 

To prove I.5, on the angles of an isosceles triangle ABC,  Euclid does not superpose the triangle on itself, 
 but relies instead on a more elaborate construction that has the advantage of also handling external angles.   
It is legitimate to interchange the two rays forming an angle at A, but it seems that the assumptions described in our analysis of I.4
are not in themselves sufficient for proving that the triangles AFC and AGB in I.5 are congruent, something essential for continuing the proof.
Here one needs to invoke the hidden assumption about angles that lie in different triangles but generate the same angle, in the modern sense.
That assumption is also used in the proofs of the two cases of I.7, where angles are decomposed into two and magnitudes of angles in different 
triangles are compared.  Although Proclus made no objections to any of the procedures just criticized, some of his comments about 
the results and their proofs are of interest.

We claim that the Fifth Common Notion, in its application to angles, is an unnecessarily strong assumption.
In the alternative approach developed below, it is a consequence of two axioms about right angles, one arguably a form of  
the Fourth Postulate that all right angles are equal.  The basic results needed for developing plane Euclidean geometry
can then be proved, albeit in a different order and, in the initial stages, via distinctly different methods.

\section{Right angles}

A well-known way to define right angles and develop their properties uses 
isosceles triangles, which are easily constructed from a chord of a circle and its center.  
By triangle transport, any such ABC with AB \cg AC can be mapped to itself, interchanging B and C.  
As the map preserves congruences on the line through BC, it fixes the midpoint M of BC.
The triangles AMB and AMC are then copies of each other.  
Whenever such triangles exist,  the angles \an AMB and \an AMC (as unions of two rays) form what we call {\em paired right angles} at M, 
on the side of the line BC containing A.  A {\em right angle} is one that is part of a pair on at least one of its sides.
A possible refinement, one we choose not to make, would use directed angles, with one of each kind in a pair of right angles.
By triangle transport, paired right angles, not necessarily unique, exist at each point along any line, on either side of that line (a location).
The purpose and intended meaning of Euclid's Fourth Postulate, that all right angles are equal,  have often been questioned (see [DR3]).
The postulate does not by itself assert the uniqueness  of right angles (as rays) at each location. We will soon suggest a satisfying interpretation.

The isosceles triangle ABC and M, used just above, can be copied to some DBC on the other side of BC, again producing paired 
right angles at M. The diagram now consists of a quadrilateral ABDC with all sides equal (a rhombus), one diagonal BC, and its midpoint M.

In the easy case that M lies on the other diagonal AD and is its midpoint, there are four angles at M, where any two consecutive ones are 
paired right angles.  The rhombus is uniquely determined by any two of its vertices on different diagonals, along with M.
These diagonals lie on lines that are said to be {\em perpendicular}, a relation that is symmetric.

It could be expected that a rhombus must have these properties in any model of geometric interest.
This, along with the fundamental initial results of the {\em Elements}, will be established after adding two axioms to our weak system. 
The first can be regarded as the Fourth Postulate in a precise form that asserts the location-independence of properties of 
individual right angles.

\medskip
{ \bf Axiom R1: \ If \an BAC and \an EDF are right angles, with
\nl \hbox to 87pt{\hfill}  AB \cg DE and AC \cg DF, then BC \cg EF.}
\medskip

Given just those right angles at A and D as pairs of rays, any triangle BAC with that angle at A can be copied to a 
triangle GDH with G on $\ora{DE}$ and H on $\ora{DF}$.
 One easy consequence of the axiom is that every right-angled triangle can be paired with another to form an isosceles triangle,
but it is not yet clear that all copies of right-angled triangles are right-angled.
As the axioms given so far may not be strong enough to establish basic results of the {\em Elements} that exclude seemingly absurd 
configurations that will soon arise, an explicit restriction is added.  A striking way to do so is as follows.

\medskip
{ \bf Axiom R2: At most one angle of a triangle is a right angle.} 
\medskip

Many other choices could be made.  One, in modern language, would require that a line and a circle cannot have infinitely many points in common.
It is not inconceivable that this could fail in some model, say one with a line that intersects a circle in an infinitesimal arc.

\begin{thm} In every rhombus the diagonals intersect at their midpoints and are mutually perpendicular.
\end{thm}

{\em Proof:} Let M be the midpoint of the diagonal BC of the rhombus ABDC. Assume for the sake of contradiction that A, M and D are not collinear.
By Axiom R1, the right triangle AMB can be transported to some D\pr MB, where D\pr\ lies on the ray $\ora{MD}$. 
By Axiom R1,  ABD\pr C is a rhombus. In its interior  BC and AD\pr\  intersect at a point P.   
The triangles AMP and D\pr MP have right angles at M and AM \cg D\pr M, so P is the midpoint  of AD\pr, by Axiom R1.
These same triangles also have right angles at P, since AMD\pr\ is isosceles.  By Axiom R2, this cannot occur.
Thus M must lie on AD, so the diagonals of the given arbitrary rhombus ABDC intersect at the midpoint M of BC.
By symmetry, M is also the midpoint of AD.  This is the easy case discussed above, so the diagonals are mutually perpendicular.  
 \eop

\begin{cor} Through each point M on a line $l$ there passes a unique line $m$ perpendicular to $l$.
\end{cor} 

{\em Proof:}  By triangle transport,  at any point M of $l$, there are paired right angles formed as usual from an isosceles triangle BAC,
 with BC along $l$, and a ray $\ora{MA}$.  On the other side of $l$,  consider a pair at M formed from another such triangle B\pr DC\pr.
After altering D so that  MD \cg MA, preserving the ray $\ora{MD}$, a rhombus ABDC is formed, by Axiom R1. 
As proven just above, the rays $\ora{MD}$ and $\ora{MA}$ must be opposites, forming a common line $m$ perpendicular to $l$.
It is then not possible for any other ray at M to form, with $l$, a pair of right angles.  
\eop

\begin{cor} Distinct circles have at most two points in common.
\end{cor} 

{\em Proof:}
This is Euclid's Prop.~III.10, whose proof using perpendicular bisectors of chords is justified by the result just above. 
A circle has a unique center, as it lies on the perpendicular bisector of every chord.
This center is then determined by any three distinct points on that circle, so no other circle can contain more than two of those points.
\eop

\section{Rigidity and angle measure}

Angles within triangles can be transported, but it is important to obtain a proof of Prop.~ I.7, a result about rigidity, within our system of axioms.
This is now effected by using, once more, configurations on each side of a line.
 
\begin{thm} If AB \cg A\pr B and  AC \cg A\pr C, with A and A\pr\ on the same side of the line through B and C, then A\pr\ = A. 
\end{thm}

{\em Proof:}  The triangle ABC can be copied to a triangle DBC (in that order)  on the other side of BC.  Then each of A, A\pr\  and D lies on 
two circles, one with center B and the other with center C.  Since those circles cannot have more than two common points,  A\pr\ must coincide with A. 
\eop
 
Conversely, Axiom R1 and the congruence-preserving construction of Theorem~1 allow existing triangles to be copied to any location,
assuming only that there is a way (by construction or axioms) to drop perpendiculars to lines and to erect them at a given point on a line.
This leads to the following significantly stronger form of rigidity for triangle transport:

\begin{thm}
If a triangle ABC is copied to A\pr B\pr C\pr, mapping points D and E on the lines formed by the sides of ABC to D\pr\ and E\pr, then D\pr E\pr\ \cg  DE.
\end{thm}

{\em Proof:}  It can be supposed that DE is not collinear with any side of ABC, and the lines through AB and through BC contain D and E, respectively.  
Let F be the foot of the perpendicular from A to the line $l$ through BC, with images $l'$ and F\pr\  under the copying process.  
By  Axiom R1 and Theorem 3, A\pr\ must coincide with the point A\st\ on the same side of $l'$ which lies
on the perpendicular to $l'$ at F\pr, with  A\st F\pr\ \cg AF.
Then A\pr E\pr\ \cg AE, since either E = F or the segments in question are hypotenuses of corresponding right triangles.  
The initial suppositions imply that  E $\neq$  B, so the same copying process can be applied to the triangle ABE and D, giving  D\pr E\pr\ \cg  DE.   
\eop

As an immediate consequence, an angle formed from two rays at a vertex A can be copied uniquely to any specified location, 
for the previous result guarantees that the process does not depend on the choice of triangle ABC having that given angle at A.
An angle can also be copied to itself, interchanging its two rays.  Just as for segments, angles and also triangles are defined to be congruent
precisely when they are copies of each other, but here no further axioms are needed to ensure that these are equivalence relations. 
It is also clear that each congruence class of angles has a unique member in a fixed reference location, say with one fixed ray at O 
along a line $l$ and the other ray on the chosen side of $l$.   Routine verifications show that the set of congruence classes of 
angles supports a natural ordering and well-defined partial operations of sum and difference, with the expected properties.  

Enough has now been established to be able to follow in sequence the Propositions in the {\em Elements} that have not already been proven.
It may be necessary to assume in places that certain circle-circle or line-circle pairs intersect.  
This does not exhaust the potential of using Euclidean axiomatics within a set-theoretic framework, as will next be shown.

One idea allows the magnitude of an angle to be defined numerically instead of as an equivalence class, so that 
magnitudes of sums and differences of angles, where defined, are the corresponding sums and differences of the magnitudes of the angles.
In other words, angles, via their magnitudes, will satisfy the first four of the five common notions. 
Since any angle can be bisected repeatedly,  among those that can be constructed are the ones whose ratio with a right angle 
is expressible as $\frac{m}{2^n}$, where $0 < m < 2^{n+1}$.  Following ancient traditions, a right angle has magnitude 90\dg, 
so the magnitudes of the above angles form a dense set in the interval from 0\dg\ to 180\dg.  
By analogy with Dedekind's construction of the reals using cuts, inspired by the theory of proportions of Eudoxus, 
every angle can be assigned a magnitude in the obvious way.
There may exist {\em infinitesimal} angles, ones that are less than all repeated bisections of right angles.  
These have magnitude 0\dg, and their supplements have magnitude 180\dg.  Every other angle has magnitude strictly between 
these extremes, and the magnitude determines the angle up to adding or subtracting an infinitesimal angle.

\section{Coordinatization and rigid motion}

We end by obtaining, with little effort,  a result about coordinate change and rigid motion in planar absolute geometry.
This involves only a superficial study of the analytic geometry of such planes, in contrast to more sophisticated approaches to
coordinatization that lead to far deeper results, as outlined in [Gr, App.~B].

Near the end of Section~4,  it was described how the common notions for segments in a geometry determine, in a concrete fashion,
a certain ordered abelian group $T$.  In addition, any ordered line $l$ and point O on $l$ determine a unique order-preserving bijection 
between $T$ and the points of $l$ (a coordinatization) such that the zero of $T$ maps to O and the magnitudes of subsegments AB of $l$,
as positive elements of $T$, can be read off as the difference of the coordinates of A and B.

By a {\em system of axes} we mean an ordered pair $(l,m)$ of perpendicular ordered lines, dividing the plane into four quadrants.
Each of these lines is coordinatized by $T$ in the way just mentioned, with 0 at the point O where the lines intersect.
Relative to any such system, the plane can be coordinatized by assigning to each point P an ordered pair of elements of $T$, 
where the first member is the coordinate of the foot of the perpendicular from P to $l$, and the second uses $m$ in place of $l$.
Points are determined uniquely by such rectilinear coordinates, as perpendiculars erected at given points of $l$ or $m$ are unique,
but a pair of elements of $T$  need not determine a point.
 
A {\em rigid motion} of the plane is a permutation of the points that preserves congruences of segments, 
hence also congruences of angles (from triangle transport). 
The image of a system of axes under any such motion is also a system of axes, where a point on one of the original axes maps to 
the point of the new axis that has the same coordinate.   To study coordinate changes and rigid motions, it is more convenient 
at this stage to use a form of polar coordinates that uses congruence classes of angles rather than angular measure, and either 
extends the range of permitted angles or notes the quadrants in which points or rays lie.  The point at the origin can safely be ignored.

\begin{thm} Any map as above between systems of axes $(l,m)$ and $(l',m')$ extends uniquely to a rigid motion of the plane.
\end{thm}
 
{ \em Proof:}  Both $l$ and $m$ are separated by O (their intersection point) into two opposite rays.  Any point P not on $l$ or $m$ lies in the interior 
of one of the right angles formed by those rays, say $\ora{OA}$ and $\ora{OB}$, and $\ora{OP}$ separates the right angle into two angles.
By triangle transport, there is a unique point P\pr\  whose polar coordinates in the system $(l',m')$ coincide with those of P in the original system.
This extension of the given map clearly defines a bijection $X \to X'$ of the plane.   To see that it is rigid, consider distinct points P and  Q.
Ignoring the easy case where O, P, Q are collinear, the polar coordinates of P and Q determine the angle $\angle POQ$.  In the new system,
P\pr\ and Q\pr\ determine the same angle, so the triangles POQ and P\pr O\pr Q\pr\ are congruent by the SAS criterion (Prop.~I.4).  
In particular, P\pr Q\pr\ \cg PQ.
 \eop

Rectilinear coordinates could have been used to define the bijective map, but 
wthout considering angles it would not be clear that this map preserves congruences.
One natural question is easily answered:

\begin{prop} The set of pairs of $T$ that correspond to points is independent of the system of axes that defines the coordinatization. 
\end{prop}

{ \em Proof:}  This sort of independence holds for systems of polar coordinates.  From the rigidity of pairs of right triangles with a common
hypotenuse, each with a side obtained by dropping a perpendicular to an axis,  one can see that in every system of axes the polar 
coordinates of a point determine the same rectilinear ones.
\eop 

The geometries for which every pair in $T \times T$ determines a point of the plane are precisely those that satisfy the  
Lotschnitt axiom, a weakening of the parallel postulate about which much is known.  See  [PS] or, for brief details, [Wi].
There are some obvious restrictions on which subsets of $T \times T$ can be the coordinate set of a geometry from which 
$T$ is defined, but we do not pursue the questions that this raises.

\pdfbookmark[1]{References}{}
\begin{sloppypar} \printbibliography \end{sloppypar}

\end{document}